\documentclass[reqno]{amsart}
\usepackage{amssymb}
\usepackage{amsmath}
\usepackage{amsfonts}
\usepackage{graphicx}
\usepackage{environ}

\newcommand{\vertiii}[1]{{\left\vert\kern-0.25ex\left\vert\kern-0.25ex\left\vert #1 
    \right\vert\kern-0.25ex\right\vert\kern-0.25ex\right\vert}}

\newtheorem{thm}{Theorem}[section]

\newtheorem{cor}[thm]{Corollary}
\newtheorem{lem}[thm]{Lemma}
\newtheorem{rem}[thm]{Remark}
\newtheorem{con}[thm]{Conjecture}

\newtheorem{note}[thm]{Note}
\newtheorem{que}[thm]{Question}

\newtheorem{ex}[thm]{Example}

\newtheorem{define}[thm]{Definition}
\numberwithin{equation}{section}

\newcommand{\CC}{\mathbb{C}}

 \newcommand{\MM}{\mathbb{M}}
 
 \newcommand{\PP}{\mathbb{P}}
\newcommand{\ie}{\text{,~i.e.,~}} 
 \newcommand{\BOX}{\hfill $\Box$}

\def\alert#1{\smallskip{\hskip\parindent\vrule\vbox{\advance\hsize-2\parindent\hrule\smallskip\parindent.4\parindent\narrower\noindent#1\smallskip\hrule}\vrule\hfill}\smallskip}

\newcommand{\mythe}[2]{
\begin{thm}
\ 
\label{#1}#2
\end{thm}}

\newcommand{\mycon}[2]{
\begin{con}
\ 
\label{#1} #2
\end{con}}

\newcommand{\mylem}[2]{
\begin{lem}
\ 
\label{#1} #2
\end{lem}}

\newcommand{\mycor}[2]{ 
\begin{cor}
\ 
\label{#1}#2
\end{cor}}

\newcommand{\myque}[2]{
\begin{que}
\ 
\label{#1}\rm #2
\end{que}}

\begin{document}

\title[Norm inequalities involving   geometric means]{Norm inequalities involving   geometric means}
\author[Freewan]{Shaima'a Freewan$^{1}$}
\address{$^{1}$Department of Mathematics, Yarmouk University, Irbid, Jordan}
\email{shyf725@gmail.com}
\author[Hayajneh]{Mostafa Hayajneh$^{2}$\ }
\address{$^{2}$Department of Mathematics, Yarmouk University, Irbid, Jordan}
\email{hayaj86@yahoo.com}


\begin{abstract}
Let $A_i$ and $B_i$ be positive definite matrices  for every $i=1,\cdots,m.$ Let $Z=[Z_{ij}]$ be the  block matrix, where  $Z_{ij}=B_i^{^\frac{1}{_2}}\left(\displaystyle\sum_{k=1}^mA_k\right)B_j^{^\frac{1}{_2}}$   for every $ i,j=~1,\cdots,m$.  It is shown that $$\left|\left|\left|\sum_{i=1}^m\left(A_i^{s}\sharp B_i^{s}\right)^r\right|\right|\right|\leq\left|\left|\left|Z^{^\frac{sr}{_2}}\right|\right|\right|\leq\left|\left|\left|\left(\left(\scalebox{0.85}{$\displaystyle\sum_{i=1}^mA_i$}\right)^\frac{srp}{_4}\left(\scalebox{0.85}{$\displaystyle\sum_{i=1}^mB_i$}\right)^\frac{srp}{_2}\left(\scalebox{0.85}{$\displaystyle\sum_{i=1}^mA_i$}\right)^\frac{srp}{_4}\right)^{\frac{1}{_p}}\right|\right|\right|,$$ for all $s\geq2$, for all $p>0$  and $r\geq1$ such that  $rp\geq1$ and for all unitarily invariant norms.

 This result generalizes the results in \cite{ONIR} and gives  an affirmative answer  to a conjecture in  \cite{OACRT}  for all $s\geq2$  and  for all $p>0$  and $r\geq1$ such that  $rp\geq1$ and $t=\frac{1}{2}$.
 This result  also  leads directly to Dinh, Ahsani, and Tam's conjecture in  \cite{GAI} and  proves Audenaert's  result in  \cite{ANIFP}.   
\end{abstract}

\maketitle

\section{Introduction}

Throughout this paper,  the symbol $\MM_n(\CC)$  refers to the set consisting of complex matrices with size $n$.  
If $\langle x, Ax\rangle \geq 0$  for each vector $x\in\CC^n$, then we can say that a Hermitian  matrix $A$ is positive semidefinite.
If $\langle x, Ax\rangle >0$  for each vector $x\in\CC^n-\{0\}$, then we can say that a Hermitian  matrix $A$ is positive definite.
   We use the symbol $\PP_n$ for the set consisting of   positive definite matrices in $\MM_n(\CC)$.
 We use the symbol  $\vertiii{.}$ for  unitarily invariant norms on the space $\MM_n(\CC)$.

For all $A,B\in\PP_n,$  the matrix $A\sharp B$ is called the   geometric mean of $A$ and $B$  given by \begin{eqnarray}\label{A471}
A\sharp B=A^{^\frac{1}{_2}}\left( A^{^\frac{-1}{_2}}BA^{^\frac{-1}{_2}}\right)^{^\frac{1}{_2}} A^{^\frac{1}{_2}}.
\end{eqnarray}
 The matrix $A\sharp_t B$ is called the $t$-geometric mean of $A$ and $B$ given by \begin{eqnarray}\label{A4732}
A\sharp_t B=A^{^\frac{1}{_2}}\left( A^{^\frac{-1}{_2}}BA^{^\frac{-1}{_2}}\right)^{t} A^{^\frac{1}{_2}},~\hspace{0.1cm}\forall ~t\in[0,1].
\end{eqnarray} 
The matrix $|A|$  is called the absolute value of $A$ defined as $|A|=(A^*A)^{^\frac{1}{_2}}.$ 
It is well known that the subadditivity relation  \begin{eqnarray}\label{E2}
f(a + b)\leq f(a) + f(b),
\end{eqnarray} holds for any $a, b \geq0$ and  for any  nonnegative concave functions $f$ on $[0,\infty).$

Bourin and Uchiyama \cite[Theorem 1.1]{AMS} discovered a noncommutative version of this inequality for all positive semidefinite matrices.  This version states that  if $A,B\in\MM_n(\CC)$ are  both positive semidefinite matrices, then for all unitarily invariant norms, we have \begin{eqnarray}\label{E1}
\left|\left|\left|f\left(A+B\right)\right|\right|\right|\leq\left|\left|\left|f\left(A\right)+f\left(B\right)\right|\right|\right|.
\end{eqnarray} 

This result prompted Bourin \cite{MSIAB} to pose the following question.
\myque{E3}{If  $A,B\in\MM_n(\CC)$  are positive semidefinite matrices, is it true that \begin{eqnarray*}
\left|\left|\left|A^{p+q}+B^{p+q}\right|\right|\right|\leq\left|\left|\left|(A^{p}+B^{p})(A^{q}+B^{q})\right|\right|\right|,~ \forall~ p,q>0?
\end{eqnarray*} }
Hayajneh and Kittaneh \cite{TIAAQ}  gave an affirmative answer to this question for the Hilbert-Schmidt and trace norms.  Audenaert \cite{ANIFP}  gave  an affirmative answer to this question by proving  \begin{eqnarray}\label{E4}
\left|\left|\left|\sum_{i=1}^m\left(A_i B_i\right)\right|\right|\right|\leq\left|\left|\left|\left(\sum_{i=1}^mA_i^\frac{1}{_2}B_i^\frac{1}{_2}\right)^2\right|\right|\right|\leq\left|\left|\left|\left(\sum_{i=1}^mA_i\right)\left(\sum_{i=1}^mB_i\right)\right|\right|\right|,
\end{eqnarray} for any $A_i,B_i\in\PP_n,~i=1,\cdots,m$ such that $A_iB_i=B_iA_i$ and  for all  unitarily invariant norms.
This result also confirms  a conjecture that Hayajneh and Kittaneh proved in \cite{TIAAQ}. Additionally, Lin \cite{ROTRR}  gave a different proof for inequality  (\ref{E4}). 

Hayajneh and Kittaneh \cite{TIAAQ2} provided  a  generalization of  inequality (\ref{E4}) by proving
\begin{eqnarray}\label{E5}
\left|\left|\left|\sum_{i=1}^m\left(A_i B_i\right)\right|\right|\right|\leq\left|\left|\left|\left(\sum_{i=1}^mA_i^\frac{1}{_2}B_i^\frac{1}{_2}\right)^2\right|\right|\right|\leq\left|\left|\left|\left(\sum_{i=1}^mA_i\right)^\frac{1}{2}\left(\sum_{i=1}^mB_i\right)\left(\sum_{i=1}^mA_i\right)^\frac{1}{2}\right|\right|\right|,\hspace{0.3cm}
\end{eqnarray}for any $A_i,B_i\in\PP_n,~i=1,\cdots,m,$ such that $A_iB_i=B_iA_i$ and  for all  unitarily invariant norms. 

Dinh, Ahsani, and Tam  \cite[Theorem~$3.1$]{GAI}  made a noncommutative version  of inequality (\ref{E5}).  According to this version, if  $A_i,B_i\in\PP_n$ for all  $i=1,\cdots,m,$ then for all unitarily invariant norms, we have 
 \begin{eqnarray}\label{A1}
\left|\left|\left|\sum_{i=1}^m\left(A_i\sharp B_i\right)^2\right|\right|\right|\leq\left|\left|\left|\left(\sum_{i=1}^mA_i\right)^\frac{1}{2}\left(\sum_{i=1}^mB_i\right)\left(\sum_{i=1}^mA_i\right)^\frac{1}{2}\right|\right|\right|.~~~
\end{eqnarray}

Dinh \cite{AIF} has generalized inequality (\ref{A1}). This generalization states that if  $A_i,B_i\in\PP_n$ for all $i=1,\cdots,m$ and if $t\in[0,1],$ $p>0$  and $r\geq1$, then
\begin{eqnarray}\label{A2}
\left|\left|\left|\sum_{i=1}^m\left(A_i\sharp_t B_i\right)^r\right|\right|\right|\leq\left|\left|\left|\left(\left(\sum_{i=1}^mA_i\right)^\frac{(1-t)rp}{2}\left(\sum_{i=1}^mB_i\right)^{trp}\left(\sum_{i=1}^mA_i\right)^\frac{(1-t)rp}{2}\right)^\frac{1}{p}\right|\right|\right|,\hspace{0.3cm}
\end{eqnarray}for all unitarily invariant norms. 

Inequality (\ref{E5}) prompted also the authors in \cite[page 787]{GAI} to  provide the following conjecture.
\mycon{CON5}{If $A_i,B_i\in\PP_n$ for all  $i=1,\cdots,m,$ then \begin{eqnarray}\label{A492}
\left|\left|\left|\sum_{i=1}^m\left(A_i^2\sharp B_i^2\right)\right|\right|\right|\leq\left|\left|\left|\left(\sum_{i=1}^mA_i\right)^\frac{1}{2}\left(\sum_{i=1}^mB_i\right)\left(\sum_{i=1}^mA_i\right)^\frac{1}{2}\right|\right|\right|,
\end{eqnarray}for  all unitarily invariant norms.}

The authors  in \cite[Corollary 3.3]{GAI}  provided an affirmative answer to Conjecture \ref{CON5}  for the trace norm.  Freewan and Hayajneh \cite{ONIR}  settled Conjecture  \ref{CON5} in its full generality by treating the following more general conjecture. 
\mycon{CON2}{
If $A_i,B_i\in\PP_n$ for all $i=1,\cdots,m,$ and  $t\in[0,1],$ then $$\left|\left|\left|\sum_{i=1}^m(A_i^2\sharp_t B_i^2)^r\right|\right|\right|\leq\left|\left|\left|\left(\left(\sum_{i=1}^mA_i\right)^{(1-t)rp}\left(\sum_{i=1}^mB_i\right)^{2trp}\left(\sum_{i=1}^mA_i\right)^{(1-t)rp}\right)^{\frac{1}{p}}\right|\right|\right|,$$ for all  $2r\geq1$ and for all $p>0$ and  for  all unitarily invariant norms.}
 Freewan and Hayajneh \cite{ONIR} gave an affirmative answer to this conjecture  for the trace norm and $t=\frac{1}{2}$. 
 Moreover, they answered  the above conjecture affirmatively for the cases
 \begin{eqnarray}\label{AAA}
 t=\frac{1}{2},~r\geq1,~ p>0 \mathrm{~such~that~}  rp\geq1
\end{eqnarray}
  and for all unitarily invariant norms. 
Conjecture \ref{CON5} is particularly answered by this result if we take $p=1$ and $r=1$ in (\ref{AAA}).

In  \cite{OACRT}, Freewan and Hayajneh  gave an alternative proof for the cases in (\ref{AAA}) without using the method of majorization.  In other words, they proved   $$\left|\left|\left|\sum_{i=1}^m(A_i^2\sharp B_i^2)^r\right|\right|\right|\leq\left|\left|\left|Z^{^r}\right|\right|\right|\leq\left|\left|\left|\left(\left(\sum_{i=1}^mA_i\right)^{\frac{rp}{_2}}\left(\sum_{i=1}^mB_i\right)^{rp}\left(\sum_{i=1}^mA_i\right)^{\frac{rp}{_2}}\right)^{\frac{1}{p}}\right|\right|\right|,$$for any $A_i,B_i\in\PP_n$, $i=1,\cdots,m,$  for any $p>0$, $r\geq1$ such that $rp\geq1$,  and for  all unitarily invariant norms and where $Z=[Z_{ij}]$ is the  block matrix such that  $Z_{ij}=B_i^{^\frac{1}{_2}}\left(\displaystyle\sum_{k=1}^mA_k\right)B_j^{^\frac{1}{_2}}$   for every $ i,j=~1,\cdots,m$.  
  
 One may ask the following natural question: What about if we replace $(A_i^2\sharp_t B_i^2)^r$ by $(A_i^s\sharp_t B_i^s)^r$ where $sr\geq1$ in Conjecture \ref{CON2}? Actually, this is the following more general conjecture posed by  Freewan and Hayajneh in  \cite{OACRT}. 
 \mycon{CON6}{Let $A_i,B_i\in\PP_n$ for all $i=1,\cdots,m$ and let $t\in[0,1]$.~Then $$\left|\left|\left|\sum_{i=1}^m\left(A_i^{s}\sharp_t B_i^{s}\right)^r\right|\right|\right|\leq\left|\left|\left|\left(\left(\displaystyle\sum_{i=1}^mA_i\right)^\frac{(1-t)srp}{_2}\left(\displaystyle\sum_{i=1}^mB_i\right)^{tsrp}\left(\displaystyle\sum_{i=1}^mA_i\right)^\frac{(1-t)srp}{_2}\right)^{\frac{1}{_p}}\right|\right|\right|,$$ for all  $p>0 $,  for all  $r>0$ and  $s>0$ such that $sr\geq1$ and for all unitarily invariant norms.}
 In this paper, we give an affirmative answer to Conjecture \ref{CON6} when  $t=\frac{1}{2},$  $s\geq2$, $p>0 $ and $r\geq1$~such~that~$rp\geq~1.$  In other words, we prove   $$\left|\left|\left|\sum_{i=1}^m\left(A_i^{s}\sharp B_i^{s}\right)^r\right|\right|\right|\leq\left|\left|\left|Z^{^\frac{sr}{_2}}\right|\right|\right|\leq\left|\left|\left|\left(\left(\displaystyle\sum_{i=1}^mA_i\right)^\frac{srp}{_4}\left(\displaystyle\sum_{i=1}^mB_i\right)^\frac{srp}{_2}\left(\displaystyle\sum_{i=1}^mA_i\right)^\frac{srp}{_4}\right)^{\frac{1}{_p}}\right|\right|\right|,$$ for any $A_i,B_i\in\PP_n$, $i=1,\cdots,m,$ for any $s\geq2$,  for any $p>0$, $r\geq1$ such that $rp\geq1$,  and for  all unitarily invariant norms  and where $Z=[Z_{ij}]$ is the  block matrix such that  $Z_{ij}=B_i^{^\frac{1}{_2}}\left(\displaystyle\sum_{k=1}^mA_k\right)B_j^{^\frac{1}{_2}}$   for every $ i,j=~1,\cdots,m$.

 \section{Preliminaries and Definitions}
This section begins with the following two lemmas, which will be used to prove our main results.
 \mylem{Lem14}{Let $A,B\in\PP_n$ and let $q\geq1$ and $p>0$.  Then $$\left|\left|\left|(BAB)^{^{pq}}\right|\right|\right|\leq \left|\left|\left|(B^{^q}A^{^q}B^{^q})^{^p}\right|\right|\right|.$$}{\it Proof.~}See \cite[page 258]{MA} and \cite[Lemma $2. 12$]{ONIR}.\BOX

 \mylem{Lem10}{If $Z=[Z_{i,j}]\in\MM_{mn}(\CC)$ is a block matrix  such that $Z_{i,j}\in\MM_n(\CC)$ is  a normal matrix for all $i,j\in\{1,2,\cdots,m\}$ or $Z$ is Hermation,  then $$|||Z|||\leq\left|\left|\left|\displaystyle\sum_{i,j=1}^m|Z_{i,j}|\right|\right|\right|,$$for every unitarily invariant norms.}{\it Proof.~}See \cite[page 7]{MSIAB}.\BOX

 Let $A,B\in\MM_n(\CC).$ Then the symbol $A\oplus B$ indicates the  direct sum  of $A$ and $B$  defined as $$A\oplus B=\left[\begin{array}{ccc}
 A & 0\\
0 & B
\end{array}\right]\in\MM_{2n}(\CC).$$ Let $\left|\left|\left|.\right|\right|\right|$ be a unitarly invariat norm on $\MM_{2n}(\CC).$ Then we define $  \left|\left|\left|A\right|\right|\right|=\left|\left|\left|A\oplus0\right|\right|\right|$.
 Let us now provide three useful lemmas to prove Theorem  \ref {The3}.  The second lemma is the matrix version of the celebrated H\"{o}lder inequality.
  \mylem{Lem1}{If $A,B\in\PP_n$,  then there exists a unitary~ $U\in~\MM_n(\CC)$ such that 
$A\sharp B = A^{^\frac{1}{_2}}UB^{^\frac{1}{_2}}.$}
{\it Proof.~}See \cite[page 108]{PDM}.\BOX

 \mylem{Lem3}{If $A,B\in\MM_n(\CC)$, then for all unitarily invariant norms, we have$$|||XY||| \leq|||~|X|^q~|||^{^\frac{1}{q}}|||~|Y|^s~|||^{^\frac{1}{s}},$$ for all $q>1$ and  $s>1$ such  that $\frac{1}{q}+\frac{1}{s}=1$.}{\it Proof.~}See \cite[page 95]{MA}.\BOX 
 
 \mylem{Lem2}{Let $A, B\in\MM_n(\CC)$  such that the product
$AB$ is normal. Then~for every unitarily invariant norms, we have $|||AB||| \leq|||BA|||.$}{\it Proof.~}See \cite[page 253]{MA}.\BOX

The following lemmas are needed to prove  Theorem \ref {The1}.

 \mylem{Lem6}{Let $A,B\in\PP_n$ and $r\geq1.$  Then for all unitarily invariant norms $$\left|\left|\left|A\right|\right|\right|\leq\left|\left|\left|B\right|\right|\right|\Longrightarrow\left|\left|\left|A^{^r}\right|\right|\right|\leq \left|\left|\left|B^{^r}\right|\right|\right|.$$}
{\it Proof.~}See \cite[Lemma $2. 9$]{ONIR}.\BOX

\mylem{Lem12}{Let $Y\in\MM_n(\CC)$ and let $a\geq0.$ Then $$|||(Y^*Y)^a||| =|||(YY^*)^a|||, $$ for all unitarily invariant norms.}
{\it Proof.~}Let $Y\in\MM_n(\CC)$ and let $a\geq0.$ Let $Y=UP$ be the polar decomposition of $Y$, where  $U$ is unitary and $P=|Y|.$  
Now, note that $$(YY^*)^a=(UPPU^*)^a=UP^{2a}U^*=U(PU^*UP)^aU^*=U(Y^*Y)^aU^*.$$
Thus $|||(Y^*Y)^a||| =|||(YY^*)^a|||.$
  This completes the proof.\BOX

\mylem{Lem15}{Let $A_i\geq 0$  for all $i=1,\cdots,m.$ Let $f:[0,\infty)\longrightarrow[0,\infty)$ be a convex function with $f(0)=0$.
Then for all unitarily invariant norms, we have
 \begin{eqnarray*}
 \left|\left|\left|\sum_{i=1}^mf(A_i)\right|\right|\right|\leq\left|\left|\left|f\left(\sum_{i=1}^mA_i\right)\right|\right|\right|.
\end{eqnarray*}}
{\it Proof.~}See \cite[Theorem 1.2]{AMS}.\BOX

The following two lemmas are taken from  \cite[Theorem 1]{OACRT}. We choose to add their proofs so that this paper becomes self-contained.
 \mylem{Lem4}{Let $A_i,B_i\in\PP_n$ for all $i=1,\cdots,m$.  Let $Z=[Z_{ij}]$ be the  block matrix such that $Z_{ij}=B_i^{^\frac{1}{_2}}\left(\displaystyle\sum_{k=1}^mA_k\right)B_j^{^\frac{1}{_2}}$ for all $i,j=1,\cdots,m$ and let \begin{eqnarray*}
Y=\left[\scalebox{0.9}{$\begin{array}{ccccc}
B_1^{^\frac{1}{_2}}&0&\cdots&0\\
\vdots&\vdots&\ddots&\vdots\\
B_m^{^\frac{1}{_2}}&0&\cdots&0
\end{array}$}\right]\left[\scalebox{0.9}{$\begin{array}{ccc}
A_1^{^\frac{1}{_2}}&\cdots&A_m^{^\frac{1}{_2}}\\
0&\cdots&0\\
\vdots&\ddots&\vdots\\
0&\cdots&0
\end{array}$}\right].
\end{eqnarray*}   Then  $Z=YY^*$.} {\it Proof.~}Note that  \begin{eqnarray}\label{eqn110}
YY^*=\left[\scalebox{0.85}{$\begin{array}{ccccc}
B_1^{^\frac{1}{_2}}&0&\cdots&0\\
\vdots&\vdots&\ddots&\vdots\\
B_m^{^\frac{1}{_2}}&0&\cdots&0
\end{array}$}\right]\left[\begin{array}{ccc}
\displaystyle\sum_{k=1}^mA_k&0\\
0&0
\end{array}\right]\left[\scalebox{0.85}{$\begin{array}{ccc}
B_1^{^\frac{1}{_2}}&\cdots&B_m^{^\frac{1}{_2}}\\
0&\cdots&0\\
\vdots&\ddots&\vdots\\
0&\cdots&0
\end{array}$}\right]=Z.
\end{eqnarray}Thus $YY^*=Z.$
 This completes the proof.\BOX

Two matrices $ A, B \in  \MM_{n} $  are unitarily equivalence if there is a unitary matrix $ U \in \MM_{n} $ such that $ B=UAU^* $, and we write $ A \cong B $.

\mylem{Lem5}{Let $A_i,B_i\in\PP_n$ for all $i=1,\cdots,m$ and let $Z=[Z_{ij}]$ be the  block matrix such that $Z_{ij}=B_i^{^\frac{1}{_2}}\left(\displaystyle\sum_{k=1}^mA_k\right)B_j^{^\frac{1}{_2}}$ for all $i,j=1,\cdots,m$.  Then $$Z \cong \left( \left(\sum_{i=1}^mA_i\right)^{\frac{1}{_2}}\left(\sum_{i=1}^mB_i\right)\left(\sum_{i=1}^mA_i\right)^{\frac{1}{_2}}\right) \oplus 0.$$}
 {\it Proof.~}Let $$K=\left[\begin{array}{ccccc}
B_1^{^\frac{1}{_2}}&0&\cdots&0\\
\vdots&\vdots&\ddots&\vdots\\
B_m^{^\frac{1}{_2}}&0&\cdots&0
\end{array}\right]\left[\begin{array}{ccc}
\left(\displaystyle\sum_{i=1}^mA_i\right)^\frac{1}{_2}&0\\
0&0
\end{array}\right].$$ Then $$Z=KK^*\cong K^*K= \left( \left(\sum_{i=1}^mA_i\right)^{\frac{1}{_2}}\left(\sum_{i=1}^mB_i\right)\left(\sum_{i=1}^mA_i\right)^{\frac{1}{_2}}\right) \oplus 0.$$
This completes the proof.\BOX

 \section{Main Results}\label{Sec3}
In this section, we introduce our main results. Lemma  \ref {Lem10} and Lemma  \ref {Lem13} play an   important role in proving~Lemma~\ref{The2}.  In the proof of the following lemma, we need to use  Lemma \ref{Lem14}. 
\mylem{Lem13}{Let $A,B\in\PP_n$  and let $U\in\MM_n(\CC)$ be a unitary. Then for all unitarily invariant norms and for all $q\geq1$, we have\begin{eqnarray*}
\left|\left|\left|~|AUB|^q~\right|\right|\right|\leq\left|\left|\left|~\left|A^{q}UB^{q}\right|~\right|\right|\right|.
\end{eqnarray*}}
{\it Proof.~}Let $A,B\in\PP_n$  and let $U\in\MM_n(\CC)$ be a unitary. 
We are to prove that for all unitarily invariant norms and for all $q\geq1$, we have\begin{eqnarray*}
\left|\left|\left|~|AUB|^q~\right|\right|\right|\leq\left|\left|\left|~\left|A^{q}UB^{q}\right|~\right|\right|\right|.
\end{eqnarray*}
To see this, note that 
\begin{eqnarray*}
\left|\left|\left|~|AUB|^q~\right|\right|\right|
&=&\left|\left|\left|~(BU^*AAUB)^\frac{q}{_2}~\right|\right|\right|\\
&=&\left|\left|\left|~(BU^*A^2UB)^\frac{q}{_2}~\right|\right|\right|\\
&\leq&\left|\left|\left|~(B^q(U^*A^2U)^qB^q)^\frac{1}{_2}~\right|\right|\right|\hspace{0.3cm}\text{(by  Lemma \ref{Lem14})}\\
&=&\left|\left|\left|~(B^qU^*A^{2q}UB^q)^\frac{1}{_2}~\right|\right|\right|\\
&=&\left|\left|\left|~((A^{q}UB^q)^*(A^{q}UB^q))^\frac{1}{_2}~\right|\right|\right|\\
&=&\left|\left|\left|~|A^{q}UB^q|~\right|\right|\right|.
\end{eqnarray*}
So   for all unitarily invariant norms and for all $q\geq1$, we have$$\left|\left|\left|~|AUB|^q~\right|\right|\right|\leq\left|\left|\left|~|A^{q}UB^q|~\right|\right|\right|.$$
 This completes the proof.\BOX

 The following lemma plays an   important role in proving~Theorem~\ref{The3}.
 \mylem{The2}{Let $A_i,B_i\in\PP_n$  and let $U_i\in\MM_n(\CC)$ be a unitary   for all $i=~1,\cdots,m.$ Suppose that $s>1$ and $q>1$~such~that $\frac{1}{s}+\frac{1}{q}=1.$ If   $$X=\left[\scalebox{0.9}{$\begin{array}{ccc}
A_1^{^\frac{s-1}{_2}}U_1B_1^{^\frac{s-1}{_2}}&~&0\\
~&\ddots&~\\
0&~&A_m^{^\frac{s-1}{_2}}U_mB_m^{^\frac{s-1}{_2}}
\end{array}$}\right],$$ then for all unitarily invariant norms, we have\begin{eqnarray*}
\left|\left|\left|~|X|^q~\right|\right|\right|\leq\left|\left|\left|\displaystyle\sum_{i=1}^m\left|A_i^{^\frac{s}{_2}}U_iB_i^{^\frac{s}{_2}}\right|~~\right|\right|\right|.
\end{eqnarray*}}
{\it Proof.~} Let $A_i,B_i\in\PP_n$  and let $U_i\in\MM_n(\CC)$ be a unitary matrix for all $i=1,\cdots,m.$   Suppose that $s>1$ and $q>1$~such~that $\frac{1}{s}+\frac{1}{q}=1\ie(s-1)q=s.$ Let   $X\in\MM_{nm}(\CC)$ be given by $$X=\left[\scalebox{0.9}{$\begin{array}{ccc}
A_1^{^\frac{s-1}{_2}}U_1B_1^{^\frac{s-1}{_2}}&~&0\\
~&\ddots&~\\
0&~&A_m^{^\frac{s-1}{_2}}U_mB_m^{^\frac{s-1}{_2}}
\end{array}$}\right].$$  We are to prove that for all unitarily invariant norms, we have\begin{eqnarray*}
\left|\left|\left|~|X|^q~\right|\right|\right|\leq\left|\left|\left|\displaystyle\sum_{i=1}^m\left|A_i^{^\frac{s}{_2}}U_iB_i^{^\frac{s}{_2}}\right|~~\right|\right|\right|.
\end{eqnarray*}
 Now, note that \begin{eqnarray*}
\left|\left|\left|~|X|^q~\right|\right|\right| &=&\left|\left|\left|~\left|\left[\scalebox{0.9}{$\begin{array}{ccc}
A_1^{^\frac{s-1}{_2}}U_1B_1^{^\frac{s-1}{_2}}&~&0\\
~&\ddots&~\\
0&~&A_m^{^\frac{s-1}{_2}}U_mB_m^{^\frac{s-1}{_2}}
\end{array}$}\right]\right|^q~\right|\right|\right|\\
&=&\left|\left|\left|~\left|\left[\scalebox{0.9}{$\begin{array}{ccc}
A_1^{^\frac{s-1}{_2}}&~&0\\
~&\ddots&~\\
0&~&A_m^{^\frac{s-1}{_2}}
\end{array}$}\right]\left[\scalebox{0.9}{$\begin{array}{ccc}
U_1&~&0\\
~&\ddots&~\\
0&~&U_m
\end{array}$}\right]\left[\scalebox{0.9}{$\begin{array}{ccc}
B_1^{^\frac{s-1}{_2}}&~&0\\
~&\ddots&~\\
0&~&B_m^{^\frac{s-1}{_2}}
\end{array}$}\right]\right|^q~\right|\right|\right|\\
&\leq&\left|\left|\left|~\left|\left[\scalebox{0.9}{$\begin{array}{ccc}
A_1^{^\frac{s-1}{_2}}&~&0\\
~&\ddots&~\\
0&~&A_m^{^\frac{s-1}{_2}}
\end{array}$}\right]^q\left[\scalebox{0.9}{$\begin{array}{ccc}
U_1&~&0\\
~&\ddots&~\\
0&~&U_m
\end{array}$}\right]\left[\scalebox{0.9}{$\begin{array}{ccc}
B_1^{^\frac{s-1}{_2}}&~&0\\
~&\ddots&~\\
0&~&B_m^{^\frac{s-1}{_2}}
\end{array}$}\right]^q\right|~\right|\right|\right|\\
&&\hspace{8cm}\text{(by  Lemma \ref{Lem13})}\\
  &=&\left|\left|\left|~\left[\scalebox{0.9}{$\begin{array}{ccc}
\left|A_1^{^\frac{(s-1)q}{_2}}U_1B_1^{^\frac{(s-1)q}{_2}}\right|&~&0\\
~&\ddots&~\\
0&~&\left|A_m^{^\frac{(s-1)q}{_2}}U_mB_m^{^\frac{(s-1)q}{_2}}\right|
\end{array}$}\right]~\right|\right|\right|\\
 &=&\left|\left|\left|~\left[\scalebox{0.9}{$\begin{array}{ccc}
\left|A_1^{^\frac{s}{_2}}U_1B_1^{^\frac{s}{_2}}\right|&~&0\\
~&\ddots&~\\
0&~&\left|A_m^{^\frac{s}{_2}}U_mB_m^{^\frac{s}{_2}}\right|
\end{array}$}\right]~\right|\right|\right|\\
&\leq&\left|\left|\left|\displaystyle\sum_{i=1}^m\left|A_i^{^\frac{s}{_2}}U_iB_i^{^\frac{s}{_2}}\right|~~\right|\right|\right|.\hspace{0.3cm}\text{(by  Lemma \ref{Lem10})}
 \end{eqnarray*}
  So for all unitarily invariant norms, we have $ \left|\left|\left|~|X|^q~\right|\right|\right|\leq\left|\left|\left|\displaystyle\sum_{i=1}^m\left|A_i^{^\frac{s}{_2}}U_iB_i^{^\frac{s}{_2}}\right|~~\right|\right|\right|.$
  This completes the proof.\BOX

Theorem~\ref{The3} and Theorem~\ref{The1} are the main ingredients to get our main result, namely, Theorem~\ref{The5}.
\mythe{The3}{ Let $A_i,B_i\in\PP_n$ for all $i=1,\cdots,m$ and let $Z=[Z_{ij}]$ be the  block matrix such that $Z_{ij}=B_i^{^\frac{1}{_2}}\left(\displaystyle\sum_{k=1}^mA_k\right)B_j^{^\frac{1}{_2}}$ for all $i,j=1,\cdots,m$. Then for all unitarily invariant norms and for all $s\geq1$, we have \begin{eqnarray}
\left|\left|\left|\displaystyle\sum_{i=1}^m\left(A_i^{s}\sharp B_i^{s}\right)\right|\right|\right|\leq\left|\left|\left|Z^{^\frac{s}{2}}\right|\right|\right|.
\end{eqnarray}
  }
  {\it Proof.~} Let $A_i,B_i\in\PP_n$ for all $i=~1,\cdots,m$ and let $s>1$. Then, using  Lemma \ref{Lem1}, we have  $A_i^{s}\sharp B_i^{s}=A_i^{^\frac{s}{_2}}U_iB_i^{^\frac{s}{_2}}$ for some unitary $U_i\in\MM_n(\CC)$ and for all $i=1,\cdots,m.$  Let $X,Y\in\MM_{mn}(\CC)$ be given by \begin{eqnarray*}
X=\left[\scalebox{0.8}{$\begin{array}{ccc}
B_1^{^\frac{s-1}{_2}}U_1^*A_1^{^\frac{s-1}{_2}}&~&0\\
~&\ddots&~\\
0&~&B_m^{^\frac{s-1}{_2}}U_m^*A_m^{^\frac{s-1}{_2}}
\end{array}$}\right]
\mathrm{~and~}
Y=\left[\scalebox{0.8}{$\begin{array}{ccccc}
B_1^{^\frac{1}{_2}}&0&\cdots&0\\
\vdots&\vdots&\ddots&\vdots\\
B_m^{^\frac{1}{_2}}&0&\cdots&0
\end{array}$}\right]\left[\scalebox{0.9}{$\begin{array}{ccc}
A_1^{^\frac{1}{_2}}&\cdots&A_m^{^\frac{1}{_2}}\\
0&\cdots&0\\
\vdots&\ddots&\vdots\\
0&\cdots&0
\end{array}$}\right].
\end{eqnarray*} 
   We are to prove that $
\left|\left|\left|\displaystyle\sum_{i=1}^m\left(A_i^{s}\sharp B_i^{s}\right)\right|\right|\right|\leq~\left|\left|\left|Z^{^\frac{s}{2}}\right|\right|\right|$ for all unitarily invariant norms.
  To see this, note that
\begin{eqnarray*}
&&\left|\left|\left|\sum_{i=1}^m\left(A_i^{s}\sharp B_i^{s}\right)\right|\right|\right|\\
&=&\left|\left|\left|\sum_{i=1}^m\left(A_i^{^\frac{s}{_2}}U_iB_i^{^\frac{s}{_2}}\right)\right|\right|\right|\\
&=&\left|\left|\left|\left(\sum_{i=1}^m\left(A_i^{^\frac{s}{_2}}U_iB_i^{^\frac{s}{_2}}\right)\right)\oplus0\right|\right|\right|\hspace{0.3cm}\\
&=&\left|\left|\left|\left[\begin{array}{ccc}
\displaystyle\sum_{i=1}^m\left(A_i^{^\frac{s}{_2}}U_iB_i^{^\frac{s}{_2}}\right)&0\\
0&0
\end{array}\right]\right|\right|\right|\\
&=&\left|\left|\left|\left[\scalebox{0.9}{$\begin{array}{ccc}
A_1^{^\frac{1}{_2}}&\cdots&A_m^{^\frac{1}{_2}}\\
0&\cdots&0\\
\vdots&\ddots&\vdots\\
0&\cdots&0
\end{array}$}\right]\left[\scalebox{0.9}{$\begin{array}{ccc}
A_1^{^\frac{s-1}{_2}}U_1B_1^{^\frac{s-1}{_2}}&~&0\\
~&\ddots&~\\
0&~&A_m^{^\frac{s-1}{_2}}U_mB_m^{^\frac{s-1}{_2}}
\end{array}$}\right]\left[\scalebox{0.9}{$\begin{array}{ccccc}
B_1^{^\frac{1}{_2}}&0&\cdots&0\\
\vdots&\vdots&\ddots&\vdots\\
B_m^{^\frac{1}{_2}}&0&\cdots&0
\end{array}$}\right]\right|\right|\right|\nonumber\\
&\leq&\left|\left|\left|\left[\scalebox{0.9}{$\begin{array}{ccc}
A_1^{^\frac{s-1}{_2}}U_1B_1^{^\frac{s-1}{_2}}&~&0\\
~&\ddots&~\\
0&~&A_m^{^\frac{s-1}{_2}}U_mB_m^{^\frac{s-1}{_2}}
\end{array}$}\right]\left[\scalebox{0.9}{$\begin{array}{ccccc}
B_1^{^\frac{1}{_2}}&0&\cdots&0\\
\vdots&\vdots&\ddots&\vdots\\
B_m^{^\frac{1}{_2}}&0&\cdots&0
\end{array}$}\right]\left[\scalebox{0.9}{$\begin{array}{ccc}
A_1^{^\frac{1}{_2}}&\cdots&A_m^{^\frac{1}{_2}}\\
0&\cdots&0\\
\vdots&\ddots&\vdots\\
0&\cdots&0
\end{array}$}\right]\right|\right|\right|\nonumber\\
&&\hspace{10cm}\text{(by  Lemma \ref{Lem2})}\nonumber
\end{eqnarray*}
\begin{eqnarray*}
&=&\left|\left|\left|XY\right|\right|\right|\label{eqn632}\\
&\leq&|||~|X|^q~|||^{^\frac{1}{q}}|||~|Y|^s~|||^{^\frac{1}{s}}\hspace{0.3cm}\text{(by  (Lemma \ref{Lem3}))}\nonumber\\
&\leq&\left|\left|\left|\displaystyle\sum_{i=1}^m\left|A_i^{^\frac{s}{_2}}U_iB_i^{^\frac{s}{_2}}\right|~~\right|\right|\right|^{\frac{1}{q}}|||~|Y|^s~|||^{\frac{1}{s}}\hspace{0.3cm} \text{(by Lemma \ref{The2})}\nonumber\\
&=&\left|\left|\left|\displaystyle\sum_{i=1}^m\left(A_i^{s}\sharp B_i^{s}\right)\right|\right|\right|^{^\frac{1}{q}}|||~|Y|^s~|||^{^\frac{1}{s}}.\nonumber
\end{eqnarray*}
This shows that \begin{eqnarray*}
\left|\left|\left|\sum_{i=1}^m\left(A_i^{s}\sharp B_i^{s}\right)\right|\right|\right|\leq\left|\left|\left|\displaystyle\sum_{i=1}^m\left(A_i^{s}\sharp B_i^{s}\right)\right|\right|\right|^{^\frac{1}{q}}|||~|Y|^s~|||^{^\frac{1}{s}}.
\end{eqnarray*}
Thus \begin{eqnarray*}
\left|\left|\left|\displaystyle\sum_{i=1}^m\left(A_i^{^s}\sharp B_i^{^s}\right)\right|\right|\right|^{1-\frac{1}{q}}\leq|||~|Y|^s~|||^{^\frac{1}{s}}.
\end{eqnarray*}
 But $1-\frac{1}{q}=\frac{1}{s},$ so  for all unitarily invariant norms, we have
 \begin{eqnarray}\label{eqn6}
\left|\left|\left|\displaystyle\sum_{i=1}^m\left(A_i^{s}\sharp B_i^{s}\right)\right|\right|\right|\leq|||~|Y|^s~|||.
\end{eqnarray}
Now, for all unitarily invariant norms, note that \begin{eqnarray*}
\left|\left|\left|~|Y|^s~\right|\right|\right|\nonumber
&=&\left|\left|\left|~(Y^*Y)^{^\frac{s}{_2}}~\right|\right|\right|\nonumber\\
&=&\left|\left|\left|~(YY^*)^{^\frac{s}{_2}}~\right|\right|\right|\hspace{0.3cm}\text{(by Lemma  \ref{Lem12})}\nonumber\\
&=&\left|\left|\left|Z^{^\frac{s}{2}}\right|\right|\right|.\hspace{0.3cm}\text{(by (Lemma \ref{Lem4}))}\label{SSSSSSSSSSSS2}
\end{eqnarray*}
 So  for all unitarily invariant norms, we have\begin{eqnarray}\label{eqn7}
\left|\left|\left|~|Y|^s~\right|\right|\right|=\left|\left|\left|Z^{^\frac{s}{2}}\right|\right|\right|.
\end{eqnarray} From (\ref{eqn6}) and (\ref{eqn7}), we get  
\begin{eqnarray*}
\left|\left|\left|\displaystyle\sum_{i=1}^m\left(A_i^{s}\sharp B_i^{s}\right)\right|\right|\right|\leq\left|\left|\left|Z^{^\frac{s}{2}}\right|\right|\right|,
\end{eqnarray*}for all unitarily invariant norms. So we are done. For the case $s=1$,  using (\ref{eqn632}) and (\ref{eqn7}), we have $$\left|\left| \left| \sum_{i=1}^mA_i\sharp B_i\right|\right|\right|\leq\left|\left| \left|Y\right|\right|\right|=\left|\left|\left|Z^{^\frac{1}{2}}\right|\right|\right|,$$ for all unitarily invariant norms.  This completes the proof.\BOX

  \mythe{The1}{Let $A_i,B_i\in\PP_n$ for all $i=1,\cdots,m$.  Let $Z=[Z_{ij}]$ be the  block matrix such that $Z_{ij}=B_i^{^\frac{1}{_2}}\left(\displaystyle\sum_{k=1}^mA_k\right)B_j^{^\frac{1}{_2}}$ for all $i,j=1,\cdots,m$. Then for all unitarily invariant norms and for all $s\geq2,$ we have $$\left|\left|\left|Z^{^\frac{s}{_2}}\right|\right|\right|\leq\left|\left|\left|\left(\sum_{i=1}^mA_i\right)^{\frac{s}{_4}}\left(\sum_{i=1}^mB_i\right)^{\frac{s}{_2}}\left(\sum_{i=1}^mA_i\right)^{\frac{s}{_4}}\right|\right|\right|.$$}{\it Proof.~}Let $A_i,B_i\in\PP_n$ for all $i=1,\cdots,m$ and let $s\geq2$.  Let $Z=[Z_{ij}]$ be the  block matrix such that $Z_{ij}=B_i^{^\frac{1}{_2}}\left(\displaystyle\sum_{k=1}^mA_k\right)B_j^{^\frac{1}{_2}}$ for all $i,j=1,\cdots,m$.  
We have \begin{eqnarray*}
\left|\left|\left|Z^{^\frac{s}{2}}\right|\right|\right|&=&\left|\left|\left|\left(\left(\sum_{i=1}^mA_i\right)^{\frac{1}{_2}}\left(\sum_{i=1}^mB_i\right)\left(\sum_{i=1}^mA_i\right)^{\frac{1}{_2}}\right)^{\frac{s}{_2}}\right|\right|\right|\hspace{0.3cm}\text{(by (Lemma \ref{Lem5}))}\\
&\leq&\left|\left|\left|\left(\sum_{i=1}^mA_i\right)^{\frac{s}{_4}}\left(\sum_{i=1}^mB_i\right)^{\frac{s}{_2}}\left(\sum_{i=1}^mA_i\right)^{\frac{s}{_4}}\right|\right|\right|.\\
&&\hspace{2.5cm}\left(\text{by Lemma \ref{Lem14} and because $\frac{s}{2}\geq1$}\right)
\end{eqnarray*}
So  for all unitarily invariant norms, we have\begin{eqnarray*}
\left|\left|\left|Z^{^\frac{s}{_2}}\right|\right|\right|\leq\left|\left|\left|\left(\sum_{i=1}^mA_i\right)^{\frac{s}{_4}}\left(\sum_{i=1}^mB_i\right)^{\frac{s}{_2}}\left(\sum_{i=1}^mA_i\right)^{\frac{s}{_4}}\right|\right|\right|.
\end{eqnarray*} This completes the proof.\BOX

We combine Theorem \ref{The3}  and Theorem \ref{The1} to get the following corollary.
\mycor{cor4}{Let $A_i,B_i\in\PP_n$ for all $i=1,\cdots,m$.  Let $Z=[Z_{ij}]$ be the  block matrix such that $Z_{ij}=B_i^{^\frac{1}{_2}}\left(\displaystyle\sum_{k=1}^mA_k\right)B_j^{^\frac{1}{_2}}$ for all $i,j=1,\cdots,m$. Then for all unitarily invariant norms and for all $s\geq2,$ we have \begin{eqnarray*}
\left|\left|\left|\displaystyle\sum_{i=1}^m\left(A_i^{s}\sharp B_i^{s}\right)\right|\right|\right|\leq\left|\left|\left|Z^{^\frac{s}{_2}}\right|\right|\right|\leq\left|\left|\left|\left(\sum_{i=1}^mA_i\right)^{\frac{s}{_4}}\left(\sum_{i=1}^mB_i\right)^{\frac{s}{_2}}\left(\sum_{i=1}^mA_i\right)^{\frac{s}{_4}}\right|\right|\right|.
\end{eqnarray*}}
{\it Proof.~}Using Theorem \ref{The3}  and Theorem \ref{The1}, we are done.\BOX

Now,   let us  prove our  main theorem using Corollary \ref{cor4}.
\mythe{The5}{Let $A_i,B_i \in \PP_n$ for all $i=1,\cdots,m$ and let $p>0 $ and $r\geq1$ such that $rp\geq1$.  Let $Z=[Z_{ij}]$ be the  block matrix such that $Z_{ij}=B_i^{^\frac{1}{_2}}\left(\displaystyle\sum_{k=1}^mA_k\right)B_j^{^\frac{1}{_2}}$~for all $i,j=1,\cdots,m$.   Then  for all unitarily invariant norms and for all $s\geq2,$ we~have$$\left|\left|\left|\sum_{i=1}^m\left(A_i^{s}\sharp B_i^{s}\right)^r\right|\right|\right|\leq\left|\left|\left|Z^{^\frac{sr}{_2}}\right|\right|\right|\leq\left|\left|\left|\left(\left(\displaystyle\sum_{i=1}^mA_i\right)^\frac{srp}{_4}\left(\displaystyle\sum_{i=1}^mB_i\right)^\frac{srp}{_2}\left(\displaystyle\sum_{i=1}^mA_i\right)^\frac{srp}{_4}\right)^{\frac{1}{_p}}\right|\right|\right|.$$}
  {\it Proof.~}Let $A_i,B_i\in\PP_n$ for all $i=1,\cdots,m$ and  let $s\geq2$, $p>0 $ and $r\geq1$ such that $rp\geq1$.   Let $Z=[Z_{ij}]$ be the  block matrix such that $Z_{ij}=B_i^{^\frac{1}{_2}}\left(\displaystyle\sum_{k=1}^mA_k\right)B_j^{^\frac{1}{_2}}$~for all $i,j=1,\cdots,m$.   We are to prove that for all unitarily invariant norms, we have $$\left|\left|\left|\sum_{i=1}^m\left(A_i^{s}\sharp B_i^{s}\right)^r\right|\right|\right|\leq\left|\left|\left|Z^{^\frac{sr}{_2}}\right|\right|\right|\leq\left|\left|\left|\left(\left(\displaystyle\sum_{i=1}^mA_i\right)^\frac{srp}{_4}\left(\displaystyle\sum_{i=1}^mB_i\right)^\frac{srp}{_2}\left(\displaystyle\sum_{i=1}^mA_i\right)^\frac{srp}{_4}\right)^{\frac{1}{_p}}\right|\right|\right|.$$Using Corollary \ref{cor4} and the fact that $s\geq2$, we have  \begin{eqnarray*}
\left|\left|\left|\displaystyle\sum_{i=1}^m\left(A_i^{s}\sharp B_i^{s}\right)\right|\right|\right|\leq\left|\left|\left|Z^{^\frac{s}{_2}}\right|\right|\right|\leq\left|\left|\left|\left(\sum_{i=1}^mA_i\right)^{\frac{s}{_4}}\left(\sum_{i=1}^mB_i\right)^{\frac{s}{_2}}\left(\sum_{i=1}^mA_i\right)^{\frac{s}{_4}}\right|\right|\right|.
\end{eqnarray*}  From Lemma \ref{Lem6} and the fact that $r\geq1$, we get
\begin{eqnarray}\label{aaaaaaaad1}
\left|\left|\left|\left(\displaystyle\sum_{i=1}^m\left(A_i^{s}\sharp B_i^{s}\right)\right)^r\right|\right|\right|\leq\left|\left|\left|\left(Z^{^\frac{s}{_2}}\right)^r\right|\right|\right|
\leq\left|\left|\left|\left(\left(\sum_{i=1}^mA_i\right)^{\frac{s}{_4}}\left(\sum_{i=1}^mB_i\right)^{\frac{s}{_2}}\left(\sum_{i=1}^mA_i\right)^{\frac{s}{_4}}\right)^r\right|\right|\right|.\hspace{0.5cm}
\end{eqnarray}
 Using Lemma \ref{Lem15} and since the function $f:[0,\infty) \longrightarrow [0, \infty)$ given by $f(x) = x^r$
is convex, it follows that\begin{eqnarray}\label{aaaaaaaad2}
\left|\left|\left|\sum_{i=1}^m\left(A_i^{s}\sharp B_i^{s}\right)^r\right|\right|\right|\leq\left|\left|\left|\left(\sum_{i=1}^mA_i^{s}\sharp B_i^{s}\right)^r\right|\right|\right|.
\end{eqnarray}
Now, using (\ref{aaaaaaaad1}) and (\ref{aaaaaaaad2}), we have\begin{eqnarray*}
\left|\left|\left|\sum_{i=1}^m\left(A_i^{s}\sharp B_i^{s}\right)^r\right|\right|\right|&\leq&\left|\left|\left|\left(\sum_{i=1}^mA_i^{s}\sharp B_i^{s}\right)^r\right|\right|\right|\\
&\leq&\left|\left|\left|\left(Z^{^\frac{s}{_2}}\right)^r\right|\right|\right|\\
&\leq&\left|\left|\left|\left(\left(\sum_{i=1}^mA_i\right)^{\frac{s}{_4}}\left(\sum_{i=1}^mB_i\right)^{\frac{s}{_2}}\left(\sum_{i=1}^mA_i\right)^{\frac{s}{_4}}\right)^r\right|\right|\right|\\
&=&\left|\left|\left|\left(\left(\sum_{i=1}^mA_i\right)^{\frac{s}{_4}}\left(\sum_{i=1}^mB_i\right)^{\frac{s}{_2}}\left(\sum_{i=1}^mA_i\right)^{\frac{s}{_4}}\right)^{\frac{rp}{_p}}\right|\right|\right|\\
&\leq&\left|\left|\left|\left(\left(\displaystyle\sum_{i=1}^mA_i\right)^\frac{srp}{_4}\left(\displaystyle\sum_{i=1}^mB_i\right)^\frac{srp}{_2}\left(\displaystyle\sum_{i=1}^mA_i\right)^\frac{srp}{_4}\right)^{\frac{1}{_p}}\right|\right|\right|.\\
&&\hspace{3cm}\text{(by  Lemma \ref{Lem14} and because $pr\geq1$)}
\end{eqnarray*}
This shows that
$$\left|\left|\left|\sum_{i=1}^m\left(A_i^{s}\sharp B_i^{s}\right)^r\right|\right|\right|\leq\left|\left|\left|Z^{^\frac{sr}{_2}}\right|\right|\right|\leq\left|\left|\left|\left(\left(\displaystyle\sum_{i=1}^mA_i\right)^\frac{srp}{_4}\left(\displaystyle\sum_{i=1}^mB_i\right)^\frac{srp}{_2}\left(\displaystyle\sum_{i=1}^mA_i\right)^\frac{srp}{_4}\right)^{\frac{1}{_p}}\right|\right|\right|,$$ for every unitarily invariant norms. 
  This completes the proof.\BOX
  
Finally, with an eye on (\ref{A2}) and Theorem \ref{The5}, we state the following conjecture.

\mycon{CON8}{Let $A_i,B_i\in\PP_n$ for all $i=1,\cdots,m$. Let $1\leq s \leq 2$. Then $$\left|\left|\left|\sum_{i=1}^m\left(A_i^{s}\sharp B_i^{s}\right)^r\right|\right|\right|\leq\left|\left|\left|\left(\left(\displaystyle\sum_{i=1}^mA_i\right)^\frac{(1-t)srp}{_2}\left(\displaystyle\sum_{i=1}^mB_i\right)^{tsrp}\left(\displaystyle\sum_{i=1}^mA_i\right)^\frac{(1-t)srp}{_2}\right)^{\frac{1}{_p}}\right|\right|\right|,$$ for all  $p>0 $,  for all  $r\geq1$ such that $rp\geq1$ and for all unitarily invariant norms.}
 Note that, using (\ref{A2}) and Theorem \ref{The5}, the above conjecture is true for the cases $s=1$ and $s=2$, respectively. The cases $s \in (1,2)$ remains open.\\

 \noindent {\bf Conflict of interests}
 
 No potential conflict of interest was reported by the author(s).
 
\noindent {\bf Authors' contributions}

All authors wrote the main manuscript text and reviewed the manuscript.
 
 \noindent {\bf Funding}

Not applicable

\end{document}